\documentclass[a4paper,12pt]{article}
\usepackage[all]{xy}
\usepackage[T1]{fontenc}
\usepackage[dvips]{graphicx}
\usepackage{amsmath,amssymb,amsfonts,amsthm,latexsym,mathrsfs,textcomp,verbatim}
\xyoption{web}

\title{Atomic toposes and countable categoricity}
\author{{Olivia Caramello} \vspace{3 mm}\\ {\small DPMMS, University of Cambridge,}\\{\small Wilberforce Road, Cambridge CB3 0WB, UK}\\{\small O.Caramello@dpmms.cam.ac.uk}}
\date{\today}
\begin{document}
\bgroup           
\let\footnoterule\relax  
\maketitle
\flushleft  
\begin{abstract}
We give a model-theoretic characterization of the class of geometric theories classified by an atomic topos having enough points; in particular, we show that every complete geometric theory classified by an atomic topos is countably categorical. Some applications are also discussed.      
\end{abstract} 
\egroup 
\flushleft
\vspace{5 mm}


\def\Monthnameof#1{\ifcase#1\or
   January\or February\or March\or April\or May\or June\or
   July\or August\or September\or October\or November\or December\fi}
\def\today{\number\day~\Monthnameof\month~\number\year}

%
%
%
\def\pushright#1{{
   \parfillskip=0pt            
   \widowpenalty=10000         
   \displaywidowpenalty=10000  
   \finalhyphendemerits=0      
  %
   \leavevmode                 
   \unskip                     
   \nobreak                    
   \hfil                       
   \penalty50                  
   \hskip.2em                  
   \null                       
   \hfill                      
   {#1}                        
  %
   \par}}                      

\def\qed{\pushright{$\square$}\penalty-700 \smallskip}

\newtheorem{theorem}{Theorem}[section]

\newtheorem{proposition}[theorem]{Proposition}

\newtheorem{scholium}[theorem]{Scholium}

\newtheorem{lemma}[theorem]{Lemma}

\newtheorem{corollary}[theorem]{Corollary}

\newtheorem{conjecture}[theorem]{Conjecture}

\newenvironment{proofs}%
 {\begin{trivlist}\item[]{\bf Proof }}%
 {\qed\end{trivlist}}

  \newtheorem{rmk}[theorem]{Remark}
\newenvironment{remark}{\begin{rmk}\em}{\end{rmk}}

  \newtheorem{rmks}[theorem]{Remarks}
\newenvironment{remarks}{\begin{rmks}\em}{\end{rmks}}

  \newtheorem{defn}[theorem]{Definition}
\newenvironment{definition}{\begin{defn}\em}{\end{defn}}

  \newtheorem{eg}[theorem]{Example}
\newenvironment{example}{\begin{eg}\em}{\end{eg}}

  \newtheorem{egs}[theorem]{Examples}
\newenvironment{examples}{\begin{egs}\em}{\end{egs}}


\mathcode`\<="4268  
\mathcode`\>="5269  
\mathcode`\.="313A  
\mathchardef\semicolon="603B 
\mathchardef\gt="313E
\mathchardef\lt="313C

\newcommand{\app}
 {{\sf app}}

\newcommand{\Ass}
 {{\bf Ass}}

\newcommand{\ASS}
 {{\mathbb A}{\sf ss}}

\newcommand{\Bb}
{\mathbb}

\newcommand{\biimp}
 {\!\Leftrightarrow\!}

\newcommand{\bim}
 {\rightarrowtail\kern-1em\twoheadrightarrow}

\newcommand{\bjg}
 {\mathrel{{\dashv}\,{\vdash}}}

\newcommand{\bstp}[3]
 {\mbox{$#1\! : #2 \bim #3$}}

\newcommand{\cat}
 {\!\mbox{\t{\ }}}

\newcommand{\cinf}
 {C^{\infty}}

\newcommand{\cinfrg}
 {\cinf\hy{\bf Rng}}

\newcommand{\cocomma}[2]
 {\mbox{$(#1\!\uparrow\!#2)$}}

\newcommand{\cod}
 {{\rm cod}}

\newcommand{\comma}[2]
 {\mbox{$(#1\!\downarrow\!#2)$}}

\newcommand{\comp}
 {\circ}

\newcommand{\cons}
 {{\sf cons}}

\newcommand{\Cont}
 {{\bf Cont}}

\newcommand{\ContE}
 {{\bf Cont}_{\cal E}}

\newcommand{\ContS}
 {{\bf Cont}_{\cal S}}

\newcommand{\cover}
 {-\!\!\triangleright\,}

\newcommand{\cstp}[3]
 {\mbox{$#1\! : #2 \cover #3$}}

\newcommand{\Dec}
 {{\rm Dec}}

\newcommand{\DEC}
 {{\mathbb D}{\sf ec}}

\newcommand{\den}[1]
 {[\![#1]\!]}

\newcommand{\Desc}
 {{\bf Desc}}

\newcommand{\dom}
 {{\rm dom}}

\newcommand{\Eff}
 {{\bf Eff}}

\newcommand{\EFF}
 {{\mathbb E}{\sf ff}}

\newcommand{\empstg}
 {[\,]}

\newcommand{\epi}
 {\twoheadrightarrow}

\newcommand{\estp}[3]
 {\mbox{$#1 \! : #2 \epi #3$}}

\newcommand{\ev}
 {{\rm ev}}

\newcommand{\Ext}
 {{\rm Ext}}

\newcommand{\fr}
 {\sf}

\newcommand{\fst}
 {{\sf fst}}

\newcommand{\fun}[2]
 {\mbox{$[#1\!\to\!#2]$}}

\newcommand{\funs}[2]
 {[#1\!\to\!#2]}

\newcommand{\Gl}
 {{\bf Gl}}

\newcommand{\hash}
 {\,\#\,}

\newcommand{\hy}
 {\mbox{-}}

\newcommand{\im}
 {{\rm im}}

\newcommand{\imp}
 {\!\Rightarrow\!}

\newcommand{\Ind}[1]
 {{\rm Ind}\hy #1}

\newcommand{\iten}[1]
{\item[{\rm (#1)}]}

\newcommand{\iter}
 {{\sf iter}}

\newcommand{\Kalg}
 {K\hy{\bf Alg}}

\newcommand{\llim}
 {{\mbox{$\lower.95ex\hbox{{\rm lim}}$}\atop{\scriptstyle
{\leftarrow}}}{}}

\newcommand{\llimd}
 {\lower0.37ex\hbox{$\pile{\lim \\ {\scriptstyle
\leftarrow}}$}{}}

\newcommand{\Mf}
 {{\bf Mf}}

\newcommand{\Mod}
 {{\bf Mod}}

\newcommand{\MOD}
{{\mathbb M}{\sf od}}

\newcommand{\mono}
 {\rightarrowtail}

\newcommand{\mor}
 {{\rm mor}}

\newcommand{\mstp}[3]
 {\mbox{$#1\! : #2 \mono #3$}}

\newcommand{\Mu}
 {{\rm M}}

\newcommand{\name}[1]
 {\mbox{$\ulcorner #1 \urcorner$}}

\newcommand{\names}[1]
 {\mbox{$\ulcorner$} #1 \mbox{$\urcorner$}}

\newcommand{\nml}
 {\triangleleft}

\newcommand{\ob}
 {{\rm ob}}

\newcommand{\op}
 {^{\rm op}}

\newcommand{\pepi}
 {\rightharpoondown\kern-0.9em\rightharpoondown}

\newcommand{\pmap}
 {\rightharpoondown}

\newcommand{\Pos}
 {{\bf Pos}}

\newcommand{\prarr}
 {\rightrightarrows}

\newcommand{\princfil}[1]
 {\mbox{$\uparrow\!(#1)$}}

\newcommand{\princid}[1]
 {\mbox{$\downarrow\!(#1)$}}

\newcommand{\prstp}[3]
 {\mbox{$#1\! : #2 \prarr #3$}}

\newcommand{\pstp}[3]
 {\mbox{$#1\! : #2 \pmap #3$}}

\newcommand{\relarr}
 {\looparrowright}

\newcommand{\rlim}
 {{\mbox{$\lower.95ex\hbox{{\rm lim}}$}\atop{\scriptstyle
{\rightarrow}}}{}}

\newcommand{\rlimd}
 {\lower0.37ex\hbox{$\pile{\lim \\ {\scriptstyle
\rightarrow}}$}{}}

\newcommand{\rstp}[3]
 {\mbox{$#1\! : #2 \relarr #3$}}

\newcommand{\scn}
 {{\bf scn}}

\newcommand{\scnS}
 {{\bf scn}_{\cal S}}

\newcommand{\semid}
 {\rtimes}

\newcommand{\Sep}
 {{\bf Sep}}

\newcommand{\sep}
 {{\bf sep}}

\newcommand{\Set}
 {{\bf Set }}

\newcommand{\Sh}
 {{\bf Sh}}

\newcommand{\ShE}
 {{\bf Sh}_{\cal E}}

\newcommand{\ShS}
 {{\bf Sh}_{\cal S}}

\newcommand{\sh}
 {{\bf sh}}

\newcommand{\Simp}
 {{\bf \Delta}}

\newcommand{\snd}
 {{\sf snd}}

\newcommand{\stg}[1]
 {\vec{#1}}

\newcommand{\stp}[3]
 {\mbox{$#1\! : #2 \to #3$}}

\newcommand{\Sub}
 {{\rm Sub}}

\newcommand{\SUB}
 {{\mathbb S}{\sf ub}}

\newcommand{\tbel}
 {\prec\!\prec}

\newcommand{\tic}[2]
 {\mbox{$#1\!.\!#2$}}

\newcommand{\tp}
 {\!:}

\newcommand{\tps}
 {:}

\newcommand{\tsub}
 {\pile{\lower0.5ex\hbox{.} \\ -}}

\newcommand{\wavy}
 {\leadsto}

\newcommand{\wavydown}
 {\,{\mbox{\raise.2ex\hbox{\hbox{$\wr$}
\kern-.73em{\lower.5ex\hbox{$\scriptstyle{\vee}$}}}}}\,}

\newcommand{\wbel}
 {\lt\!\lt}

\newcommand{\wstp}[3]
 {\mbox{$#1\!: #2 \wavy #3$}}

\newcommand{\fu}[2]
{[#1,#2]}

\newcommand{\st}[2]
 {\mbox{$#1 \to #2$}}

\section{Some results on atomic toposes}
In this section we present some results on atomic toposes which are relevant to our characterization theorem in the second section.\\
Let us recall the following standard definition.
\begin{definition}
Let $\mathcal{E}$ be a topos. An object $A\in \cal E$ is said to be an atom of $\cal E$ if the only subobjects of $A$ (up to isomorphism) are the identity arrow $1_{A}:A\to A$ and the zero arrow $0_{A}:0\to A$, and they are distinct from each other. 
\end{definition}
The following proposition describes the behaviour of associated sheaf functors with respect to atoms.
\begin{proposition}\label{prop1}
Let $\cal E$ be a topos and $j$ a topology on it with associated sheaf functor $a_{j}:{\cal E}\to \sh_{j}(\cal E)$. If $A$ is an atom of $\cal E$ then $a_{j}(A)$ is an atom of $\sh_{j}(\cal E)$, provided that it is non-zero.
\end{proposition}
\begin{proofs}
Given a monomorphism $m:C\to a_{j}(A)$ in $\sh_{j}({\cal E})$, $m$ is a monomorphism also in $\cal E$ since the inclusion $i:\sh_{j}({\cal E})\hookrightarrow {\cal E}$ preserves monomorphisms (having a left adjoint). Now, denoted by $\eta$ the unit of the adjuction $a_{j}\dashv i$, consider the pullback
\[  
\xymatrix {
C' \ar[r]^{m'} \ar[d]  & A \ar[d]^{\eta_{A}} \\
C \ar[r]^{m} & a_{j}(A) }
\]\\ 
in $\cal E$.
The arrow $m'$ is a monomorphism in $\cal E$, being the pullback of a monomorphism, so, since $A$ is an atom of $\cal E$ we deduce that $m'$ is either (isomorphic to) the identity arrow on $A$ or the zero arrow $0_{A}$. Now, by applying $a_{j}$ to the pullback above we obtain a pullback in $\sh_{j}({\cal E})$ (as $a_{j}$ preserves pullbacks); but $a_{j}(\eta_{A})\cong 1_{a_{j}(A)}$, so $m\cong a_{j}(m')$ and $m$ is either (isomorphic to) the identity or the zero arrow on $a_{j}(A)$; of course, if $a_{j}(A)\ncong 0_{\sh_{j}({\cal E})}$ these two arrows are distinct from each other.   
\end{proofs}
We recall that an atomic topos is an elementary topos $\cal E$ which possesses an atomic geometric morphism ${\cal E}\to \Set$. We refer the reader to section C3.5 in \cite{El2} for a comprehensive treatment of the topic of atomic toposes. Here we limit ourselves to remarking the following facts.
\begin{proposition}\label{propn1}
Let $\cal E$ be a Grothendieck topos. Then\\
(i) $\cal E$ is atomic if and only if it has a generating set of atoms;\\
(ii) if $\{a_{i} \textrm{ | } i\in I\}$ is a generating set of atoms for $\cal E$ then the atoms of $\cal E$ are exactly the epimorphic images of the atoms in the generating set; in particular, $\cal E$ has only a set of (isomorphism classes of) atoms.
\end{proposition}
\begin{proofs}
(i) Suppose that $\cal E$ is atomic. Then all the subobject lattices in $\cal E$ are atomic Boolean algebras (cfr. p. 685 \cite{El2}) and hence every object of $\cal E$ can be written as a disjoint coproduct of atoms; on the other hand, there can be only a set of atoms (up to isomorphism) in $\cal E$, by the argument at the top of p. 690 \cite{El2}. Conversely, if $\cal E$ has a generating set of atoms then the full subcategory $\cal C$ of $\cal E$ on it satisfies the right Ore condition and ${\cal E}\cong \Sh({\cal C}, J_{at})$, where $J_{at}$ is the atomic topology on $\cal C$ (cfr. the discussion p. 689 \cite{El2}); so it is atomic (by Theorem C3.5.8 \cite{El2}).\\
(ii) This was remarked p. 690 \cite{El2}. 
\end{proofs}
As a consequence of Propositions \ref{prop1} and \ref{propn1}(i), we may deduce that any subtopos of an atomic Grothendieck topos $\cal E$ is atomic; indeed, the images of the atoms in a generating set of $\cal E$ via the corresponing associated sheaf functor clearly form a generating set for the subtopos. In fact, this property holds more generally at the elementary level (i.e. every subtopos of an atomic topos is atomic), by the following argument. Let $\cal E$ be an atomic topos; then, $\cal E$ being Boolean, every subtopos $\cal F$ of $\cal E$ is open (by Proposition A4.5.22 \cite{El2}) and hence the inclusion of $\cal F$ into $\cal E$ is an atomic morphism (by Proposition A4.5.1 \cite{El}); this implies that the geometric morphism ${\cal F}\to \Set$ is atomic, being the composite of two atomic morphisms (the inclusion ${\cal F}\hookrightarrow {\cal E}$ and the morphism ${\cal E}\to \Set$); so $\cal F$ is atomic. In terms of sites, if ${\cal E}\cong \Sh({\cal C}, J^{\cal C}_{at})$ (where $\cal C$ satisfies the right Ore condition and $J^{\cal C}_{at}$ is the atomic topology on it) then the subtoposes of it can be described as follows.\\     
Let $\cal F$ be a subtopos of $\cal E$; as we have already remarked, $\cal F$ must be open, that is of the form ${\cal E}/U\hookrightarrow {\cal E}$ for a subterminal object $U$ in $\cal E$. Now, by Remark C2.3.21 \cite{El2}, $U$ can be identified with a $J_{at}$-ideal on $\cal C$, that is with a collection of objects ${\cal C}'$ of $\cal C$ with the property that for any arrow $f:a\to b$ in $\cal C$, $a\in {\cal C}'$ if and only if $b\in {\cal C}'$. If we regard ${\cal C}'$ as a full subcategory of $\cal C$ then $\Sh({\cal C}, J^{\cal C}_{at})/U\cong \Sh({\cal C}', J^{{\cal C}'}_{at})$ (where $J^{{\cal C}'}_{at}$ is the atomic topology on ${\cal C}'$). Indeed, we may define an equivalence as follows. Given a object $G\to U$ in $\Sh({\cal C}, J^{\cal C}_{at})/U$, for every $c\in {\cal C}$ which does not belong to ${\cal C}'$ we must have $G(c)=\emptyset$, since we have an arrow $G(c)\to U(c)$ and $U(c)=\emptyset$; so we associate to it the restriction $G|_{{\cal C}'}$, which is a $J^{{\cal C}'}_{at}$-sheaf since $J^{{\cal C}'}_{at}$ clearly coincides with the Grothendieck topology induced by $J^{\cal C}_{at}$ on ${\cal C}'$. It is now clear that this assigment defines an equivalence between our two categories. So we have proved that the subtoposes of $\Sh({\cal C}, J^{\cal C}_{at})$ are exactly those of the form $\Sh({\cal C}', J^{{\cal C}'}_{at})$ where ${\cal C}'$ is a full subcategory of $\cal C$ with the property that for any arrow $f:a\to b$ in $\cal C$, $a\in {\cal C}'$ if and only if $b\in {\cal C}'$. Also, since the assigment sending a subterminal object in $\cal E$ to the corresponding open subtopos of $\cal E$ is a lattice isomorphism from $\Sub_{\cal E}(1)$ to the lattice of open subtoposes of $\cal E$, two such subtoposes of $\Sh({\cal C}, J^{\cal C}_{at})$ are equivalent if and only if the corresponding categories are equal (as subcategories of $\cal C$).\\ 

Next, let us consider a general category $\cal C$. We know that, provided that $\cal C$ satisfies the right Ore condition, one can define the atomic topology on $\cal C$ as the topology having as covering sieves exactly the non-empty ones. Such a topology does not exist on a general category $\cal C$ but, by analogy with it, we may define the atomic topology $J^{\cal C}_{at}$ on $\cal C$ as the smallest Grothendieck topology on $\cal C$ such that all the non-empty sieves are covering; of course, this definition specializes to the well-known one in the case $\cal C$ satisfies the right Ore condition. As stated in following proposition, the corresponding category of sheaves is an atomic topos.
\begin{proposition}\label{propat}
Let $\cal C$ be a category and $J^{\cal C}_{at}$ the atomic topology on it. Then $\Sh({\cal C}, J^{\cal C}_{at})$ is an atomic topos. 
\end{proposition}
\begin{proofs}
Let ${\cal C}'$ be the full subcategory of $\cal C$ on the objects which are not $J^{\cal C}_{at}$-covered by the empty sieve. Then, by the Comparison Lemma, we have that $\Sh({\cal C}, J^{\cal C}_{at})\cong \Sh({\cal C}', J^{\cal C}_{at}|_{{\cal C}'})$. We now prove that ${\cal C}'$ satisfies the right Ore condition and $J^{\cal C}_{at}|_{{\cal C}'}=J^{{\cal C}'}_{at}$, that is for every sieve $R$ in ${\cal C}'$, $R\neq \emptyset$ if and only if $R$ is $J^{\cal C}_{at}|_{{\cal C}'}$-covering; from this our thesis will clearly follow. In one direction, suppose that $R\neq \emptyset$. Then the sieve $\overline{R}$ generated by $R$ in $\cal C$ is obviously non-empty and, ${\cal C}'$ being a full subcategory of $\cal C$, we have that $\overline{R}\cap arr({\cal C}')=R$; so $R$ is $J^{\cal C}_{at}|_{{\cal C}'}$-covering by definition of induced topology. Conversely, suppose that $R$ is a $J^{\cal C}_{at}|_{{\cal C}'}$-covering sieve on an object $c\in {\cal C}'$. Then there exists a $J^{\cal C}_{at}$-covering sieve $H$ on $c$ in $\cal C$ such that $H\cap arr({\cal C}')=R$. Suppose $R$ be empty; then for every arrow $f$ in $H$ we have $\emptyset\in J^{\cal C}_{at}(dom(f))$. But $H$ is $J^{\cal C}_{at}$-covering so from the transitivity axiom for Grothendieck topologies it follows that $\emptyset \in J^{\cal C}_{at}(c)$, contradiction since $c\in {\cal C}'$. So we conclude that $R$ is non-empty, as required.       
\end{proofs}   
\begin{rmk}
\emph{By the transitivity axiom for Grothendieck topologies, the subcategory ${\cal C}'$ in the proof of the proposition above satisfies the property that for any arrow $f:a\to b$ in $\cal C$, $a\in {\cal C}'$ if and only if $b\in {\cal C}'$; in other words, ${\cal C}'$ is a union of connected components of $\cal C$. In particular, if ${\cal C}'\neq {\cal C}$ (i.e. $\cal C$ does not satisfy the right Ore condition) and $\cal C$ is connected then ${\cal C}'=\emptyset$, that is the topos $\Sh({\cal C}, J^{\cal C}_{at})$ is trivial.}
\end{rmk}
The following result generalizes the proposition above.
\begin{proposition}\label{propel}
Let $\cal E$ be a Grothendieck topos with a generating set $\cal L$ and $j$ be an elementary topology on $\cal E$ such that all the monomorphisms $a\to b$ in $\cal E$ where $a\ncong 0$ and $b\in {\cal L}$ are $j$-dense. Then $\sh_{j}({\cal E})$ is an atomic topos.
\end{proposition}
\begin{proofs}
By Proposition \ref{prop1}, it is enough to prove that the images of the objects of $\cal L$ via the associated sheaf functor $a_{j}$ form a generating set of objects of $\sh_{j}({\cal E})$ which are either zero or atoms. Our argument follows the lines of the proof of Proposition \ref{prop1}. Given an object $b\in {\cal L}$ and a monomorphism $m:a\to a_{j}(b)$ in $\sh_{j}({\cal E})$, consider the pullback  
\[  
\xymatrix {
a' \ar[r]^{m'} \ar[d]  & b \ar[d]^{\eta_{b}} \\
a \ar[r]^{m} & a_{j}(b) }
\]\\ 
in $\cal E$.
The arrow $m'$ is a monomorphism in $\cal E$, being the pullback of a monomorphism, so, if $a'\ncong 0$ then $m'$ is $j$-dense by our hypotheses, that is $a_{j}(m')$ is an isomorphism. But $a_{j}$ preserves pullbacks, from which it follows that $m$ is an isomorphism. If instead $a'\cong 0$ then $a\cong a_{j}(a')\cong a_{j}(0)=0_{\sh_{j}({\cal E})}$ so $m$ is the zero arrow on $a_{j}(b)$. 
\end{proofs}
\begin{rmk}
\emph{We note that Proposition \ref{propat} is the particular case of Proposition \ref{propel} when $\cal E$ is a presheaf topos $[{\cal C}^{\textrm{op}}, \Set]$, $\cal L$ is the collection of all the representables on $\cal C$ and $j$ is the elementary topology on $[{\cal C}^{\textrm{op}}, \Set]$ corresponding to the atomic topology on $\cal C$; indeed, the sieves in $\cal C$ on an object $c\in {\cal C}$ can be identified with the subobjects in $[{\cal C}^{\textrm{op}}, \Set]$ of the representable ${\cal C}(-,c)$.}
\end{rmk}
Now, let us briefly consider another approach for obtaining an atomic topos starting from a general one, based on the consideration of the atoms of the given topos.
\begin{proposition}
Let $\cal E$ be a Grothendieck topos and $\cal L$ a collection of atoms of $\cal E$, regarded as a full subcategory of $\cal E$. Then, if $J^{\cal E}_{can}$ is the canonical topology on $\cal E$, the topos $\Sh({\cal L}, J^{\cal E}_{can}|_{\cal L})$ is atomic.
\end{proposition}
\begin{proofs}
Obviously, since every arrow in $\cal L$ is an epimorphism in $\cal E$, we have $J^{\cal L}_{at}\subseteq J^{\cal E}_{can}|_{\cal L}$ so $\Sh({\cal L}, J^{\cal E}_{can}|_{\cal L})$ is a subtopos of the topos $\Sh({\cal L}, J^{\cal L}_{at})$. But $\Sh({\cal L}, J^{\cal L}_{at})$ is atomic by Proposition \ref{propat}, hence $\Sh({\cal L}, J^{\cal L}_{at})$ is atomic by the discussion following the proof of Proposition \ref{propn1}.  
\end{proofs}    
Let us now characterize the atoms of the topos $\Sh({\cal C}, J^{\cal C}_{at})$, where $\cal C$ is a category satisfying the right Ore condition.
\begin{proposition}\label{loccon}
Let $\Sh({\cal C}, J)$ be a locally connected topos, and $a_{J}:[{\cal C}^{\textrm{op}}, \Set]\to \Sh({\cal C}, J)$ be the associated sheaf functor. Then all the functors $a_{J}({\cal C}(-,c))$ are connected objects of $\Sh({\cal C}, J)$ if and only if all the constant functors ${\cal C}^{\textrm{op}} \to \Set$ are $J$-sheaves. 
\end{proposition}
\begin{proofs}
Consider the diagram
\[  
\xymatrix {
\Sh({\cal C}, J) \ar_{p}[dr] \ar^{i}[rr]  &  & [{\cal C}^{\textrm{op}}, \Set] \ar^{q}[dl]  \\
& \Set &}
\]\\
of geometric morphisms in the 2-category of Grothendieck toposes, where $p$ and $q$ are the unique geometric morphisms respectively from $\Sh({\cal C}, J)$ and $[{\cal C}^{\textrm{op}}, \Set]$ to $\Set$. Both these geometric morphisms are essential, that is their inverse image functors have left adjoints, which we indicate respectively by $p_{!}$ and $q_{!}$; indeed, $p$ is essential because by hypothesis $\Sh({\cal C}, J)$ is locally connected, while $q$ is essential by Example A4.1.4 \cite{El}. It is well-known that the representables in $[{\cal C}^{\textrm{op}}, \Set]$ are all indecomposable, so $q_{!}({\cal C}(-,c))=1$ for each $c\in {\cal C}$. Now, the condition that all the constant functors ${\cal C}^{\textrm{op}} \to \Set$ are $J$-sheaves is clearly equivalent to demanding that $q^{\ast}=i\circ p^{\ast}$ where $i$ is the inclusion $\Sh({\cal C}, J)\hookrightarrow [{\cal C}^{\textrm{op}}, \Set]$ or, passing to the left adjoints, that $q_{!}=p_{!}\circ a$ (of course, the equalities here are intended to be isomorphisms); but, since all these functors preserve colimits (having right adjoints) and every functor in $[{\cal C}^{\textrm{op}}, \Set]$ is a colimit of representables, the equality above holds if and only if $1=q_{!}({\cal C}(-,c))=p_{!}(a_{J}({\cal C}(-,c))$, that is if and only if the $a_{J}({\cal C}(-,c))$ are all connected objects of $\Sh({\cal C}, J)$.     

\end{proofs}
\begin{rmk}\label{locconrmk}
\emph{We note that for a general Grothendieck site $({\cal C}, J)$, the constant functor $\Delta{\emptyset}:{\cal C}^{\textrm{op}} \to \Set$ is a $J$-sheaf if and only if every $J$-covering sieve is non-empty, and all the constant functors $\Delta{L}:{\cal C}^{\textrm{op}} \to \Set$ for a non-empty set $L\in \Set$ are $J$-sheaves if and only if for each object $c\in {\cal C}$, all the $J$-covering sieves on $c$ are empty or connected as full subcategories of ${\cal C}/c$; in particular, the conjuction of these two conditions implies, by Theorem C3.3.10 \cite{El2}, that the topos $\Sh({\cal C}, J)$ is locally connected.}   
\end{rmk}
As a consequence of Proposition \ref{loccon} and Remark \ref{locconrmk}, we deduce that if $\cal C$ is a category satisfying the right Ore condition and $J$ is a Grothendieck topology on $\cal C$ such that every $J$-covering sieve is non-empty, then all the functors $a({\cal C}(-,c))$ are connected objects of the locally connected topos $\Sh({\cal C}, J)$. In particular, if $J^{\cal C}_{at}$ is the atomic topology on $\cal C$ then the $a({\cal C}(-,c))$ are all atoms of the atomic topos $\Sh({\cal C}, J^{\cal C}_{at})$ (since in an atomic topos the atoms are precisely the connected objects, cfr. p. 685 \cite{El2}); since they also form a generating set for the topos $\Sh({\cal C}, J^{\cal C}_{at})$, we deduce from Proposition \ref{propn1}(ii) that the atoms of $\Sh({\cal C}, J^{\cal C}_{at})$ are exactly the epimorphic images of the functors of the form $a({\cal C}(-,c))$. By using Yoneda's lemma, one can easily rephrase this condition as follows: a $J^{\cal C}_{at}$-sheaf $F$ is an atom of $\Sh({\cal C}, J^{\cal C}_{at})$ if and only if there exists an object $c\in {\cal C}$ and an element $x\in F(c)$ with the property that every natural transformation $\alpha$ from $F$ to any $J^{\cal C}_{at}$-sheaf $G$ is uniquely determined by its value $\alpha(c)(x)$ at $x$.      
\section{The characterization theorem}
In this section we prove our main characterization result concerning the geometric theories classified by an atomic topos with enough points.\\
Let us first introduce the relevant definitions and establish some basic facts. For the general background we refer the reader to \cite{El2}.\\ Concerning notation, for convenience signatures are supposed to be one-sorted throughout the whole section, but all the arguments can be easily adapted to the general many-sorted case.
 
\begin{definition}
Let $\mathbb{T}$ be a geometric theory. $\mathbb{T}$ is said to be atomic if its classifying topos $\Set[\mathbb T]$ is an atomic topos.  
\end{definition}
\begin{definition}
Let $\mathbb{T}$ be a geometric theory over a signature $\Sigma$. $\mathbb T$ is said to have enough models if for every geometric sequent $\sigma$ over $\Sigma$, $M\vDash \sigma$ for all the $\mathbb T$-models $M$ in $\Set$ implies that $\sigma$ is provable in $\mathbb T$.  
\end{definition}
Note that since the soundness theorem for geometric logic always holds (see for example Proposition D1.3.2 p. 832 \cite{El2}), the class of theories with enough models is exactly the class of geometric theories for which `the' completeness theorem holds.  

\begin{proposition}\label{enough}
Let $\mathbb{T}$ be a geometric theory over a signature $\Sigma$. Then $\mathbb T$ has enough models if and only if its classifying topos $\Set[\mathbb T]$ has enough points. 
\end{proposition}
\begin{proofs}
By definition, $\Set[\mathbb T]$ has enough points if and only if the inverse image functors $f^{\ast}$ of the geometric morphisms $f:\Set \rightarrow \Set[\mathbb T]$ are jointly conservative. Now, since the geometric morphism $f_{M}:\Set \rightarrow \Set[{\mathbb T}]$ corresponding to a $\mathbb T$-model $M$ in $\Set$ satisfies $f^{\ast}(M_{\mathbb{T}})=M$ (where $M_{\mathbb{T}}$ is the universal model of $\mathbb T$ lying in $\Set[\mathbb T]$) then it follows from Lemma D1.2.13 p. 825 \cite{El2} that if a geometric sequent $\sigma$ over $\Sigma$ is satisfied in every $\mathbb T$-model $M$ in $\Set$ then $\sigma$ is satisfied in $M_{\mathbb{T}}$, equivalently it is provable in $\mathbb T$.\\
Conversely, suppose that $\mathbb T$ has enough models. Then it is easily seen, by using an argument analogous to that employed in the proof of Proposition D3.3.13 p. 915 \cite{El2}, that $\Set[\mathbb T]$ has enough points.     
\end{proofs}
\begin{definition}
Let $\mathbb{T}$ be a geometric theory over a signature $\Sigma$. $\mathbb T$ is said to be complete if every geometric sentence $\phi$ over $\Sigma$ is $\mathbb T$-provably equivalent to $\top$ or $\bot$, but not both.  
\end{definition}
\begin{rmk}\label{rmk2}
\emph{From the topos-theoretic point of view, a geometric theory is complete if and only if its classifying topos is two-valued (to see this, it suffices to consider the syntactic representation for the classifying topos as the category of sheaves on the geometric syntactic category of the theory with respect to the `syntactic topology' on it); moreover, if $\mathbb{T}$ is atomic then its classifying topos is two-valued if and only if it is (atomic and) connected (cfr. the proof of Theorem 2.5. \cite{OC2}).}\\ 
\end{rmk}
Given a geometric theory $\mathbb T$ over a signature $\Sigma$, from now on we will denote the relation of $\mathbb T$-provable equivalence of geometric formulas over $\Sigma$ in the same context by $\stackrel{\mathbb T}{\sim}$.
\begin{definition}
Let $\mathbb{T}$ be a geometric theory over a signature $\Sigma$. $\mathbb T$ is said to be Boolean if it classifying topos is a Boolean topos. 
\end{definition}
\begin{rmk}\label{rmk1}
\emph{We recall from \cite{OC3} that a geometric theory $\mathbb T$ over a signature $\Sigma$ is a Boolean if and only if for every geometric formula $\phi(\vec{x})$ over $\Sigma$ there exists a geometric formula $\psi(\vec{x})$ over $\Sigma$ in the same context, denoted $\neg \phi(\vec{x})$, such that $\phi(\vec{x}) \wedge \psi(\vec{x})\stackrel{\mathbb T}{\sim}\bot$ and $\phi(\vec{x}) \vee \psi(\vec{x})\stackrel{\mathbb T}{\sim}\top$.\\
From this criterion, it follows that if $\mathbb T$ is a Boolean then every infinitarily disjunctive first-order formula over $\Sigma$ (i.e. an infinitary first-order formula over $\Sigma$ which do not contain infinitary conjunctions) is $\mathbb T$-provably equivalent using classical logic to a geometric formula in the same context; indeed, this can be proved by an inductive argument as in the proof of Theorem D3.4.6 p. 921 \cite{El2}.} 
\end{rmk}

\begin{definition}
Let $\mathbb{T}$ be a geometric theory over a signature $\Sigma$. Two $\mathbb T$-models (in $\Set$) $M$ and $N$ are said to be geometrically equivalent if and only if for every geometric sentence $\phi$ over $\Sigma$, $M\vDash \phi$ if and only if $N\vDash \phi$. 
\end{definition}
Let us recall that a model $M$ of a geometric theory $\mathbb T$ over a signature $\Sigma$ is said to be conservative if $M\vDash \sigma$ for every geometric sequent $\sigma$ over $\Sigma$ implies $\sigma$ provable in $\mathbb T$.\\ 
The following result represents the geometric analogue of the well-known characterization of completeness of a first-order theory in model theory.
Below, by a trivial geometric theory we mean a geometric theory in which $\bot$ is provable.
\begin{proposition}\label{prop}
Let $\mathbb{T}$ be a non-trivial Boolean geometric theory with enough models. Then the following are equivalent:\\
(i) $\mathbb T$ is complete;\\
(ii) for every geometric sentence $\phi$, either $\phi\stackrel{\mathbb T}{\sim}\top$ or $\neg \phi \stackrel{\mathbb T}{\sim} \top$;\\   
(iii) every two $\mathbb T$-models in $\Set$ are geometrically equivalent;\\
(iv) every $\mathbb T$-model $M$ in $\Set$ is conservative.\\
\end{proposition}

\begin{proofs}
(i) $\biimp$ (ii) is obvious.\\
(i) $\imp$ (iii) For any geometric sentence $\phi$ over $\Sigma$, either $\phi \stackrel{\mathbb T}{\sim} \top$, and hence $M\vDash \phi$ for all the $\mathbb T$-models, or $\phi \stackrel{\mathbb T}{\sim} \bot$, and hence $M\nvDash \phi$ for all $\mathbb T$-models; so (iii) immediately follows.\\
(iii) $\imp$ (i) Given a geometric sentence $\phi$ over $\Sigma$, since $\mathbb T$ has enough models, if $\phi \stackrel{\mathbb T}{\nsim} \top$ then there exists a $\mathbb T$-model $M$ in $\Set$ such that $\phi$ does not hold in $M$; then $\phi$ does not hold in any $\mathbb T$-model in $\Set$, these models being all geometrically equivalent. This precisely means that the geometric sequent $\phi \: \vdash_{[]} \bot$ holds in every $\mathbb T$-model in $\Set$, that is, $\mathbb T$ having enough models, $\phi \stackrel{\mathbb T}{\sim}\bot$.\\  
(iii) $\imp$ (iv) Given a geometric sequent $\phi \: \vdash_{\vec{x}} \psi$ over $\Sigma$, it is clear that for any $\mathbb T$-model $M$, $\phi \: \vdash_{\vec{x}} \psi$ holds in $M$ if and only if the infinitarily disjunctive first-order sentence $\forall \vec{x}(\phi \to \psi)$ holds in $M$. But, by Remark \ref{rmk1}, this formula is $\mathbb T$-provably equivalent using classical logic to a geometric sentence; so we conclude that if a geometric sequent is satisfied in a $\mathbb T$-model $M$ then it is satisfied in every $\mathbb T$-model in $\Set$ and hence, $\mathbb T$ having enough models, it is provable in $\mathbb T$.\\
(iv) $\imp$ (iii) is obvious.\\    
\end{proofs}
\begin{rmks}
\emph{(a) As it is clear from the proof, the equivalence (i) $\biimp$ (iii) in the proposition above holds in general for any geometric theory with enough models.\\
(b) Since every Boolean topos having enough points is atomic (Corollary C3.5.2 p. 685 \cite{El2}), the implication (i) $\imp$ (iv) in the proposition above can be seen, in view of Remark \ref{rmk2}, as the logical version of the topos-theoretic fact that every point of a connected atomic topos is a surjection (cfr. Proposition C3.5.6(ii) \cite{El2}).}   
\end{rmks}
     
\begin{definition}
Let $\mathbb{T}$ be a geometric theory over a signature $\Sigma$. A type-in-context (or, more briefly, a type) of $\mathbb T$ is any set of geometric formulas over $\Sigma$ in the same context of the form $\{\phi(\vec{x}) \textrm{ | } M\vDash \phi(\vec{a})\}$, where $M$ is a model of $\mathbb T$ in $\Set$ and $\vec{a}$ is a tuple of elements of (the underlying set of) $M$; the type $\{\phi(\vec{x}) \textrm{ | } M\vDash \phi(\vec{a})\}$ will be denoted by $S^{\mathbb{T}}_{(M,\vec{a})}$.\\
A type of $\mathbb T$ is said to be complete if it is maximal (with respect to the inclusion) in the set of all types of $\mathbb T$.\\
A type $S$ of $\mathbb T$ is said to be principal if there exists a formula $\phi(\vec{x})\in S$ such that for any geometric formula $\psi(\vec{x})$ over $\Sigma$ in the same context, $\phi(\vec{x})$ $\mathbb T$-provably implies $\psi(\vec{x})$ if (and only if) $\psi(\vec{x})\in S$; the formula $\phi(\vec{x})$ is said to be a generator of the type $S$.       
\end{definition}
 
\begin{rmk}\label{rmkcomplete}
\emph{Note that, by Proposition \ref{prop}, the notion of complete geometric theory introduced above rewrites in terms of types as follows: a non-trivial geometric theory $\mathbb T$ having enough models is complete if and only if for any two $\mathbb T$-models $M$ and $N$ in $\Set$, $S^{\mathbb{T}}_{(M,[])}=S^{\mathbb{T}}_{(N,[])}$.} 
\end{rmk}

\begin{definition}
Let $\Sigma$ be a signature, $M$ a $\Sigma$-structure and $N$ a substructure of $M$. Then $N$ is said to be a geometric substructure of $M$ if, for every geometric formula $\phi(\vec{x})$ over $\Sigma$ and any tuple of elements $\vec{a}$ (of the same length as $\vec{x}$) from $N$, $M \vDash \phi(\vec{a})$ if and only if $N \vDash \phi(\vec{a})$; equivalently, $S^{\emptyset}_{(M, \vec{a})}=S^{\emptyset}_{(N, \vec{a})}$ for any tuple $\vec{a}$ of elements of $N$ (where $\emptyset$ denotes the empty geometric theory over $\Sigma$). 
\end{definition}
\begin{rmk}\label{rmk3}
\emph{It is easy to prove by induction on the structure of geometric formulas that every geometric formula is equivalent in geometric logic to an infinitary disjunction of geometric formulas which do not contain infinitary disjunctions; since these latter formulas are in particular first-order, we may deduce that if $N$ is an elementary substructure of $M$ then $N$ is a geometric substructure of $M$; moreover, given a geometric sequent $\phi(\vec{x}) \vdash_{\vec{x}} \psi(\vec{x})$, if this sequent holds in $M$ then it also holds in $N$. Indeed, for every tuple $\vec{a}$ of elements in $N$ (of the same length as $\vec{x}$), $N\vDash \phi(\vec{a})$ implies $M\vDash \phi(\vec{a})$, which in turn implies $M\vDash \psi(\vec{a})$ and hence $N\vDash \psi(\vec{a})$ (where the first and third implications follow from the fact that $N$ is a geometric substructure of $M$). We note that this remark justifies the use of the downward L\"owenheim-Skolem theorem in the context of geometric logic; more precisely, given a geometric theory $\mathbb T$ over a signature $\Sigma$ of cardinality $|\Sigma|$, if $\mathbb T$ has a model $M$ such that $|M|\geq |\Sigma|$ then $\mathbb T$ has a model of cardinality $|\Sigma|$.}     
\end{rmk}
Below by `countable' we mean either finite or denumerable.
\begin{definition}
Let $\mathbb{T}$ be a geometric theory. Then $\mathbb T$ is said to be countably categorical if any two models of $\mathbb T$ in $\Set$ of countable cardinality are isomorphic.
\end{definition}
We remark that, by our definition, any geometric theory having no models in $\Set$ is (vacously) countably categorical.\\

The following definition is the geometric equivalent of the notion of atomic model in classical model theory. 
\begin{definition}
Let $\mathbb{T}$ be a geometric theory over a signature $\Sigma$. A model $M$ of $\mathbb T$ in $\Set$ is said to be atomic if for any tuple of elements $\vec{a}$ of $M$, the type $S^{\mathbb{T}}_{(M,\vec{a})}$ is principal and complete. 
\end{definition}
Let us recall from \cite{OC3} that a geometric theory over a signature $\Sigma$ is Boolean if and only if every geometric formula $\phi(\vec{x})$ over $\Sigma$ which is stably consistent with respect to $\mathbb T$ (i.e. such that $\phi(\vec{x})\wedge \psi(\vec{x})\stackrel{\mathbb T}{\nsim}\bot$ for every geometric formula $\psi(\vec{x})$ over $\Sigma$ in the same context) is provable in $\mathbb T$; let us also recall from \cite{El2} that a geometric theory $\mathbb T$ is atomic if and only if all the subobject lattices in the geometric syntactic category ${\cal C}_{\mathbb T}$ of $\mathbb T$ are atomic Boolean algebras (this also follows from the results in the first section by using the fact that every subobject in the classifying topos $\Set[{\mathbb T}]$ of $\mathbb T$ of an object in ${\cal C}_{\mathbb T}$ lies in ${\cal C}_{\mathbb T}$). We will make use of these characterizations in the proof of the theorem below.\\

\begin{theorem}\label{teofond}
Let $\mathbb T$ be a complete geometric theory having a model in $\Set$. Then the following are equivalent:\\
(i) $\mathbb T$ is countably categorical and Boolean\\
(ii) $\mathbb T$ is atomic\\
(iii) every $\mathbb T$-model in $\Set$ is atomic
\end{theorem}

\begin{proofs}
(i) $\imp$ (ii) By Proposition \ref{prop}, any Boolean complete geometric theory with a model in $\Set$ has enough models; so the thesis follows from the fact that every Boolean topos with enough points is atomic (Corollary C3.5.2 p. 685 \cite{El2}).\\

(ii) $\imp$ (iii) Let $M$ be a $\mathbb T$-model in $\Set$ and $\vec{a}$ be a tuple of elements of $M$; we want to prove that  $S^{\mathbb{T}}_{(M,\vec{a})}$ is principal and complete. Consider the subobject lattice $\Sub_{{\cal C}_{\mathbb T}}(\{\vec{x}.\top\})$ in the geometric syntactic category ${\cal C}_{\mathbb T}$ of $\mathbb T$, where $\vec{x}$ is a set of variables of the same length as $\vec{a}$. Since $\Sub_{_{{\cal C}_{\mathbb T}}}(\{\vec{x}.\top\})$ is an atomic Boolean algebra, we can write $\{\vec{x}.\top\}$ as a disjuction of atoms of $\Sub_{{\cal C}_{\mathbb T}}(\{\vec{x}.\top\})$; so, since $\{\vec{x}.\top\}$ obviously belongs to $S^{\mathbb{T}}_{(M,\vec{a})}$, there exists exactly one atom of $\Sub_{_{{\cal C}_{\mathbb T}}}(\{\vec{x}.\top\})$ (up to $\mathbb T$-provable equivalence) which belongs to $S^{\mathbb{T}}_{(M,\vec{a})}$; then it is clear that this atom generates the type $S^{\mathbb{T}}_{(M,\vec{a})}$. So we have proved that all the types of $\mathbb T$ are principal; it remains to verify that they are also complete. To this end, let us first observe that $\mathbb T$ is Boolean (since every atomic topos is Boolean). So, given an inclusion $S^{\mathbb{T}}_{(M,\vec{a})}\subseteq S^{\mathbb{T}}_{(N,\vec{b})}$ of types of $\mathbb T$, this inclusion must be an equality because if there were a formula $\phi(\vec{x})\in S^{\mathbb{T}}_{(N,\vec{b})}\setminus S^{\mathbb{T}}_{(M,\vec{a})}$ then, by definition of $\neg \phi(\vec{x})$, we would have $\neg \phi(\vec{x})\in S^{\mathbb{T}}_{(M,\vec{a})}$ and hence $\neg \phi(\vec{x})\in S^{\mathbb{T}}_{(N,\vec{b})}$, a contradiction.\\ 

(iii) $\imp$ (ii) Let us first prove that $\mathbb T$ is Boolean, that is every formula $\phi(\vec{x})$ which is stably consistent with respect to $\mathbb T$ is provable in $\mathbb T$. Given a $\mathbb T$-model $M$ and a tuple $\vec{a}$ of elements of $M$ of the same length as $\vec{x}$, let $\psi_{(M, \vec{a})}$ be a generator of the type $S^{\mathbb{T}}_{(M,\vec{a})}$. As we have already observed, under our hypotheses $\mathbb T$ has enough models so, since $\phi(\vec{x})\wedge \psi_{(M, \vec{a})}\stackrel{\mathbb T}{\nsim}\bot$, there exists a $\mathbb T$-model $N$ and a tuple $\vec{b}$ of elements of it (of the same length as $\vec{x}$) such that $\phi(\vec{x})$ and $\psi_{(M, \vec{a})}$ both belong to $S^{\mathbb{T}}_{(N,\vec{b})}$. Now, since $\psi_{(M, \vec{a})}$ generates the type $S^{\mathbb{T}}_{(M,\vec{a})}$, it follows that $S^{\mathbb{T}}_{(M,\vec{a})}\subseteq S^{\mathbb{T}}_{(N,\vec{b})}$ and hence, since all the types of $\mathbb T$ are complete,  $S^{\mathbb{T}}_{(M,\vec{a})}=S^{\mathbb{T}}_{(N,\vec{b})}$. This in turn implies that $\phi(\vec{x})\in S^{\mathbb{T}}_{(M,\vec{a})}$, that is $M \vDash \phi(\vec{a})$. Since the $\mathbb{T}$-model $M$ and the tuple $\vec{a}$ are arbitrary, we conclude, again by invoking the fact that $\mathbb T$ has enough models, that $\phi(\vec{x})$ is provable in $\mathbb{T}$, as required. Now that we have proved that $\mathbb T$ is Boolean, to show that $\mathbb T$ is atomic, it remains to verify that all the Boolean subobject lattices in the geometric syntactic category ${\cal C}_{\mathbb T}$ of $\mathbb T$ are atomic, equivalently for every formula $\phi(\vec{x})\stackrel{\mathbb T}{\nsim} \bot$ there exists an atom below it in the Boolean algebra $\Sub_{{\cal C}_{\mathbb T}}(\{\vec{x}.\top\})$. If $\phi(\vec{x})\stackrel{\mathbb T}{\nsim} \bot$ then, since $\mathbb T$ has enough models, there exists a $\mathbb T$-model $M$ and a tuple $\vec{a}$ of elements of it (of the same length as $\vec{x}$) such that $\phi(\vec{x})\in S^{\mathbb{T}}_{(M,\vec{a})}$. It is now enough to check that the generator $\psi_{(M, \vec{a})}$ of the type $S^{\mathbb{T}}_{(M,\vec{a})}$ is an atom of $\Sub_{{\cal C}_{\mathbb T}}(\{\vec{x}.\top\})$; this follows similarly as above by using the fact that $\mathbb{T}$ has enough models and the types of $\mathbb T$ are complete.\\

(ii) $\imp$ (i) Being atomic, $\mathbb{T}$ is Boolean, as every atomic topos is Boolean. To prove that $\mathbb T$ is countably categorical, let us distinguish two cases: either $\mathbb T$ has a finite model in $\Set$ or all the models of $\mathbb{T}$ are infinite.\\  
Let us suppose that all the models of $\mathbb{T}$ are infinite. We have to prove that any two denumerable models of $\mathbb T$ are isomorphic. We will construct explicitly such an isomorphism as in the proof of Theorem 7.2.2 p. 336 \cite{Hodges}. Let $M$ and $N$ be two models of $\mathbb{T}$ of cardinality $\aleph_{0}$. Then, $\mathbb{T}$ being complete, we have $S^{\mathbb{T}}_{(M,[])}=S^{\mathbb{T}}_{(N,[])}$ by Remark \ref{rmkcomplete}. Let us first prove by induction on $k \in \mathbb{N}$ the following fact: given tuples $\vec{a}$ and $\vec{b}$ of length $k$ respectively in $M$ and $N$ such that $S^{\mathbb{T}}_{(M, \vec{a})}=S^{\mathbb{T}}_{(N, \vec{b})}$, and an element $d\in N$ there exists an element $c\in M$ such that $S^{\mathbb{T}}_{(M, \vec{a},c)}=S^{\mathbb{T}}_{(N, \vec{b},d)}$ (and, symmetrically, given an element $c\in M$ there exists an element $d\in N$ such that $S^{\mathbb{T}}_{(M, \vec{a},c)}=S^{\mathbb{T}}_{(N, \vec{b},d)}$). Consider the type $S^{\mathbb{T}}_{(N, \vec{b},d)}$; this is principal, by our hypotheses (having already proved the implication (ii) $\imp$ (iii) in the theorem), so it is generated by a formula $\psi(\vec{x}, y)$. Now, $N \vDash (\exists y \psi(\vec{x}, y))(\vec{b})$ so since $S^{\mathbb{T}}_{(M, \vec{a})}=S^{\mathbb{T}}_{(N, \vec{b})}$ we deduce that there exists $c\in M$ such that $M \vDash \psi(\vec{a}, c)$; but $\psi(\vec{x}, y)$ is a generator of $S^{\mathbb{T}}_{(N, \vec{b},d)}$ and all the types of $\mathbb T$ are complete by our hypothesis, so we conclude that $S^{\mathbb{T}}_{(M, \vec{a},c)}=S^{\mathbb{T}}_{(N, \vec{b},d)}$, as required.       
Now, since $M$ and $N$ are geometrically equivalent by Proposition \ref{prop}, an obvious back-and-forth argument yields two sequences $(m_{0}, m_{1}, \ldots, m_{k}, ... )$ and $(n_{0}, n_{1}, \ldots, n_{k}, ... )$ enumerating respectively $M$ and $N$, such that for each $k\in \mathbb{N}$ $S^{\mathbb{T}}_{(M, m_{0}, m_{1}, \ldots, m_{k})}=S^{\mathbb{T}}_{(N, n_{0}, n_{1}, \ldots, n_{k})}$; then the map $f:M\to N$ sending each $m_{k}$ to $n_{k}$ is an isomorphism of $\mathbb{T}$-models, as it is a bijection preserving the interpretation of all the atomic formulas.\\
Let us instead suppose that $\mathbb{T}$ has a finite model $M$ in $\Set$ of cardinality $n$. Consider the geometric sequents (over $\Sigma$)\\
$\top \:\vdash_{[]}\: \exists x_{1}\ldots \exists x_{n} ( \mathbin{\mathop{\textrm{\huge $\wedge$}}\limits_{1\leq i\lt j\leq n}}x_{i}\neq x_{j})$ and\\
$ \mathbin{\mathop{\textrm{\huge $\wedge$}}\limits_{1\leq i\lt j\leq n}}x_{i}\neq x_{j} \:\vdash_{x_{1},\ldots, x_{n}, y}\: \mathbin{\mathop{\textrm{\huge $\vee$}}\limits_{1\leq i\leq n}}y=x_{i}$,\\
where for each $i$ and $j$, the expression $x_{i}\neq x_{j}$ denotes the complement of the formula $x_{i}\neq x_{j}$ in the subobject lattice $\Sub_{{\cal C}_{\mathbb{T}}}(\{x_{i}, x_{j}. \top\})$ of the geometric syntactic category ${\cal C}_{\mathbb{T}}$ of $\mathbb{T}$ (recall that, since the classifying topos of $\mathbb T$ is Boolean, these sublattices are all Boolean algebras).\\
Clearly, a model $N$ of $\mathbb{T}$ satisfies these sequents if and only if it has cardinality $n$; so in particular $M$ satisfies them. But, $\mathbb{T}$ being Boolean and complete, $M$ is a conservative model of $\mathbb T$ by Proposition \ref{prop}, so these sequents are provable in $\mathbb T$. From this, it follows that all the models of $\mathbb{T}$ have cardinality $n$. Since they are all atomic (by the implication (ii) $\imp$ (iii) in the theorem), a back-and-forth argument as above yields an isomorphism between any two models of $\mathbb{T}$.\\                            
\end{proofs}

\begin{rmks}
\emph{
(a) The equivalence (i) $\biimp$ (ii) in the theorem above generalizes the analogous result for coherent theories obtained by A. R. Blass and A. \v{S}\v{c}edrov in \cite{blasce}.\\
(b) As it is clear from the proof of the theorem above, the equivalence (ii) $\biimp$ (iii) holds in general for any geometric theory with enough models, while the implication (ii) $\imp$ (i) holds for any complete geometric theory.
}
\end{rmks}

Given a geometric theory $\mathbb{T}$ over a signature $\Sigma$, by a `quotient' of $\mathbb{T}$ we mean a geometric theory $\mathbb{T}'$ over $\Sigma$ such that every axiom of $\mathbb{T}$ is provable in $\mathbb{T}'$; if $\mathbb{T}'$ is complete, then we say that $\mathbb{T}'$ is a completion of $\mathbb{T}$.\\
Let us now describe the completions of an atomic theory $\mathbb{T}$. Since $\Sub_{{\cal C}_{\mathbb{T}}}(\{[].\top\})$ is an atomic Boolean algebra, we can write $\top$ as a disjunction $\mathbin{\mathop{\textrm{\huge $\vee$}}\limits_{i\in I}}\phi_{i}$ of geometric sentences which are atoms of $\Sub_{{\cal C}_{\mathbb{T}}}(\{[].\top\})$. Then the completions of $\mathbb T$ are precisely the theories $\mathbb{T}_{i}$ obtained from $\mathbb T$ by adding to it an axiom of the form $\top \vdash_{[]} \phi_{i}$. Indeed, by our results in the first section, a subtopos ${\cal E}\slash U$ of an atomic topos $\cal E$ is two-valued if and only if $U$ is an atom of $\cal E$; also, if $\cal E$ is atomic then we have a decomposition of $1_{\cal E}$ as a disjoint sum of atoms $\mathbin{\mathop{\textrm{\huge $\vee$}}\limits_{i\in I}}U_{i}$ of $\Sub_{\cal E}(1_{\cal E})$ and hence $\cal E$ clearly decomposes as the coproduct of the toposes ${\cal E}\slash U_{i}$ for $i\in I$. Now, if $\cal E$ is the classifying topos $\Set[\mathbb{T}]$ of an atomic theory $\mathbb T$, then the toposes appearing in such decomposition can be clearly identified as the classifying toposes $\Set[\mathbb{T}_{i}]\simeq \Set[\mathbb{T}]\slash [[\phi_{i}]]_{G}$ of the $\mathbb{T}_{i}$, where $G$ is the universal model of $\Set[\mathbb{T}]$; so we may conclude by Remark \ref{rmk2} that the completions of $\mathbb T$ are precisely the $\mathbb{T}_{i}$, and in particular that they are all atomic theories. In passing, we note that if $\cal E$ is the category $\Sh({\cal C}, J^{\cal C}_{at})$ of sheaves on a category $\cal C$ with the respect to the atomic topology $J^{\cal C}_{at}$ on it (cfr. the first section of this paper for the definition of the atomic topology on a general category), this decomposition coincides (by the results in the first section) with the decomposition of $\Sh({\cal C}, J^{\cal C}_{at})$ as the coproduct of the toposes $\Sh({\cal C}', J^{{\cal C}'}_{at})$ as ${\cal C}'$ ranges in the set of connected components of $\cal C$.\\  
By combining this discussion with Theorem \ref{teofond} we thus obtain the following result: all the completions of an atomic geometric theory are countably categorical.\\
Finally, let us indicate how it is possible to deduce from Theorem \ref{teofond} a representation result for connected atomic toposes with a point. From the proof of the theorem, it is clear that, provided that it exists, the unique (up to isomorphism) countable model $M$ of an atomic complete theory $\mathbb T$ over $\Sigma$ satisfies the following property: any two tuples from $M$ safisfy exactly the same geometric formulas over $\Sigma$ if and only if there exists an automorphism of $M$ which sends one to another. Then one can prove, by arguments analogous to those employed in the proof of Theorem 3.2 \cite{blasce}, that the classifying topos for $\mathbb T$ is equivalent to the topos of continuous $G$-sets where $G$ is the group of automorphisms of $M$ equipped with the `topology of pointwise convergence' (i.e. the topology defined by declaring a basis of neighbourhoods of the identity to consist of the subgroups $G_{\vec{a}}=\{\alpha\in G \textrm{ | $\alpha$ fixes each element of $\vec{a}$}\}$, for finite tuples $\vec{a}$ in $M$.\\           

\section{Applications}

\begin{theorem}\label{appl}
Let $\mathbb T$ be a geometric theory having a model in $\Set$ in which every stably consistent formula with respect to $\mathbb T$ is satisfied. Then $\mathbb T$ has a quotient which is complete, countably categorical, and has a model in $\Set$.
\end{theorem}
\begin{proofs}
Consider the Booleanization $\mathbb{T}'$ of the theory $\mathbb{T}$ (as it was defined in \cite{OC3}). $\mathbb{T}'$ is a geometric theory over $\Sigma$, and our hypotheses say precisely that $\mathbb{T}'$ has a model $M$ in $\Set$. Then, the geometric theory $Th(M)$ over $\Sigma$ having as axioms all the geometric sequents over $\Sigma$ which are satisfied in $M$, is complete and contains (in the obvious sense) the theory $\mathbb{T}'$; so its classifying topos $\Set[Th(M)]$ is a subtopos of the Boolean topos $\Set[{\mathbb{T}}']$, and hence it is a Boolean topos (by Proposition A4.5.22 \cite{El}). But the theory $Th(M)$ has enough models ($M$ being a conservative model for it), so $\Set[Th(M)]$ has enough points (by Proposition \ref{enough}) and hence it is atomic, by Corollary C3.5.2 \cite{El2}. Our thesis now follows from Theorem \ref{teofond}.      
\end{proofs}
\begin{rmk}\label{rmk4}
\emph{We note that if the signature of the theory $\mathbb T$ in Theorem \ref{appl} is countable then the quotient of $\mathbb T$ in the statement of the theorem has exactly one countable model in $\Set$ up to isomorphism; indeed, this follows from the downward L\"owenheim-Skolem theorem (cfr. Remark \ref{rmk3}).}
\end{rmk}

The terminology in the following result is taken from \cite{OC2}.
\begin{theorem}\label{teofraisse}
Let $\mathbb{T}$ be a theory of preshaf type such that the category $(\textrm{f.p.} {\mathbb T}\textrm{-mod}(\Set))$ satisfies the amalgamation and joint embedding properties. Then any two countable homogeneous $\mathbb T$-models in $\Set$ are isomorphic.
\end{theorem}
\begin{proofs}
As it is remarked in \cite{OC3}, the Booleanization $\mathbb{T}'$ of $\mathbb{T}$ axiomatizes the homogeneous $\mathbb T$-models. Now, we have already observed that an atomic geometric theory is complete if and only if its classifying topos is (atomic and) connected (cfr. Remark \ref{rmk2}). So $\mathbb{T}'$ is complete, since its classifying topos $\mathbb{T}'\simeq \Sh((\textrm{f.p.} {\mathbb T}\textrm{-mod}(\Set))^{\textrm{op}}, J_{at})$ is atomic and connected, by Theorems 2.5. and 2.6. in \cite{OC2}. Our thesis now follows from Theorem \ref{teofond}.
\end{proofs}

\vspace{10 mm}
{\bf Acknowledgements:} I am grateful to my Ph.D. supervisor Peter Johnstone for many useful discussions.\\

\end{document}